# Stochastic Boundedness of State Trajectories of Stable LTI Systems in the Presence of Nonvanishing Stochastic Perturbation

Peyman Azodi*, Alireza Khayatian**, Elham Jamalinia***, Peyman Setoodeh****

*Department of Electrical Engineering and Computer science, Shiraz University, p_azodi@shirazu.ac.ir

**Department of Electrical Engineering and Computer science, Shiraz University, khayatia@shirazu.ac.ir

*** Department of Electrical Engineering and Computer science, Shiraz University, e-jamalinia@shirazu.ac.ir

**** Department of Electrical Engineering and Computer science, Shiraz University, setoodeh@ieee.org

This paper studies stochastic boundedness of trajectories of a nonvanishing stochastically perturbed stable LTI system. First, two definitions on stochastic boundedness of stochastic processes are presented, then the boundedness is analyzed via Lyapunov theory. In this proposed theorem, it is shown that under a condition on the Lipchitz constant of the perturbation kernel, the trajectories remain stochastically bounded in the sense of the proposed definitions and the bounds are calculated. Also, the limiting behavior of the trajectories have been studied. At the end an illustrative example is presented, which shows the effectiveness of the proposed theory.

## I. INTRODUCTION

Perturbed systems form a large group of systems with a wide range of applications. One of the most useful classes of these systems is the set of additively perturbed systems. In the case of vanishing perturbations[1], the stability of the original system (intrinsic system) may not be affected by the perturbation. On the other hand, in the case of nonvanishing perturbations, the stability of the original system gets completely damaged. To be more precise, the trajectories started from the origin will not remain there. In this case, the most useful attempt is to analyze the boundedness of the trajectories. Once the designer knows that the states do not reach the origin, the question is "How close they will get to the origin?". In practice, nonvanishing perturbations are more frequent. This in turn, leads to the distinction of the perturbed system from the original one at the origin.

Stochastically perturbed systems are also frequent. One of the most important applications of stochastically perturbed systems is to model the influence of Gaussian measurement of a system. This model also presents the effect of weak measurement of the quantum systems in the Schrodinger picture of quantum mechanics. These features motivated the authors to analyze the stochastic boundedness of stable LTI systems perturbed by nonvanishing (nonlinear) perturbations.

---

[1] One classification of functions $g(X,t):\mathbb{R}^n \times \mathbb{R}^+ \mapsto \mathbb{R}$, takes account of whether it vanishes at the origin or not; thus, two classes of functions due to this classification are:
1. Vanishing Perturbation Kernel, $g(0,t) \equiv 0$.
2. Nonvanishing Perturbation Kernel, $\exists I \neq \phi \subseteq T, \forall t \in I : g(0,t) \neq 0$.

Valuable research on the boundedness of determined perturbed systems exists in the literature [1]. Also, the stability of stochastic systems using Lyapunov theory has been investigated in [2-6]. In [5] and [7, 8], boundedness of differential equations, based on Lyapunov approach, was analyzed. Despite their importance, there is no specific research on evaluating the stochastic bound on the trajectories of nonvanishing perturbed systems, which is the main contribution of this paper. In [9], the expectationally bounded stochastic processes have been defined and the results on the expectationally boundedness of the trajectories of a linear LTI system with stochastic nonlinear perturbation has been presented by the authors.

In this paper, first, two definitions on stochastic boundedness of stochastic processes are provided. Afterward a novel theorem is presented, which is the main contribution of this paper. By this theorem, it is shown that the LTI stable systems at the origin, which are perturbed by nonvanishing stochastic perturbations are stochastically bounded under a mild condition on the perturbation kernel (defined in subsection II.B). By this result, designers are ensured that the state trajectory is stochastically bounded; hence, they can answer the question: "How close they will get to the origin?"

In section II, the problem under consideration is formulated. First in subsection 0A, some preliminaries and definitions are presented. This subsection is followed by subsection B, in which stochastically perturbed LTI systems are described. Afterward in subsection C, two definitions on stochastically bounded stochastic processes are presented. In section III, the main results of this paper are presented in Theorem 1 and Remark 1. In section IV, an Illustrative example is presented and the effectiveness of the proposed theory is shown. The paper concludes in section V.

## II. PRELIMINARY AND PROBLEM FORMULATION

In this section, three introductory subsections are presented, which include the preliminaries and the problem formulation. Also, an introductory subsection on stochastic boundedness is presented.

### A. Preliminary and definitions

This subsection provides the required background material:

Consider a stable matrix, $A \in \mathbb{R}^{n \times n}$, in the rest of this article, $P$, which is called the *Lyapunov Complement*, represents the unique symmetric solution of the Lyapunov equation [1, 10]:

$$A^\mathrm{T} P + PA = -Q \qquad (1)$$

for any positive symmetric $Q$. In the special case that $Q = I$, this equation is called the *trivial Lyapunov equation*, and the Lyapunov complement is denoted by $\hat{P}$.

Also, $\sigma(A)$ represents the pure point spectrum of $A$, $\underline{\lambda}_A$ and $\bar{\lambda}_A$ respectively denote the eigenvalues of $A$ with the smallest and the largest absolute values:

$$\bar{\lambda}_A \triangleq \sup_i \left( |\lambda_i| \,|\, \lambda_i \in \sigma(A) \right), \quad \underline{\lambda}_A \triangleq \inf_i \left( |\lambda_i| \,|\, \lambda_i \in \sigma(A) \right) \tag{2}$$

Throughout this paper, $\|.\|$ denotes the well-known $\ell^2$ norm and $E(.)$ stands for the expectation value of a random variable.

## B. Stochastically Perturbed LTI Systems

Consider the class of nonlinearly stochastically perturbed stable LTI systems. A general form of these systems is given by the following state equation in the Ito form [3, 6]:

$$\begin{aligned} dX &= AX\,dt + g(X,t)\,dW \\ X(0) &= X_0 \quad \text{a.s} \end{aligned} \tag{3}$$

where $dW$ is an $m$ dimensional standard Wiener process defined on the underlying complete probability space $(\Omega, \Im, P)$. In this case, the stochastic process $X(t, \omega)$, which will be abbreviated by $X(t)$ or $X$, is the systems state and the strong solution of (3). This stochastic process is defined on the probability space $(\Omega, \Im, P)$ and the measurable space $(\mathbb{R}^n, B(\mathbb{R}^n))$, where $B(.)$ denotes the Borel $\sigma$–algebra of the underlying Hilbert space equipped by the Euclidean metric [11, 12]. If the time set is defined to be $T \triangleq \mathbb{R}^+ = [0, \infty)$ and $B^+ \triangleq B(T)$, then the map $X : (T \times \Omega, B^+ \otimes \Im) \to (\mathbb{R}^n, B(\mathbb{R}^n))$ is the systems response. The perturbation part of the dynamic $g(X,t)dW$ also consists of the time dependent $n \times m$ matrix valued $\Im$-measurable function $g(X,t)$, which is called Perturbation Kernel (PK), the elements of which are time-dependent real functions $g_{ij}(X,t): \mathbb{R}^n \times \mathbb{R}^+ \mapsto \mathbb{R}$. Also, it is assumed that $X_0$ is independent of the $\sigma$–algebra generated by the Wiener process $W(t)$, $\Im(W(t), 0 \leq t)$, and also $E(X_0) < \infty$.

Throughout this paper, it is assumed that the following conditions are satisfied:

**Assumption 1:**
1. $A$ is stable.
2. The PK in (3) is Lipchitz continuous on $\mathbb{R}^n$ i.e.

$$\begin{aligned} &\exists 0 \leq \gamma; \\ &\forall x, y \in \mathbb{R}^n, t \in T : \|g(x,t) - g(y,t)\| \leq \gamma \|x - y\| \end{aligned} \tag{4}$$

3. The PK is finite at the origin, i.e.

$$\sup_{t \in \mathrm{T}} \|g(0,t)\| < \infty \tag{5}$$

4. The PK admits the growth condition on $\mathbb{R}^n$ i.e.

$$\exists 0 \leq \gamma_g; \tag{6}$$
$$\forall x \in \mathbb{R}^n, t \in \mathrm{T}: \|g(x,t)\|^2 \leq \gamma_g \left(1 + \|x\|^2\right) \quad \blacksquare$$

Under the proposed assumption, the existence and uniqueness of the strong solution of (3), is guaranteed [11, 13]. Also, the proposed Assumptions will be used in the boundedness analysis. Stochastically perturbed LTI systems form a wide class of stochastic systems. In the rest of this article, boundedness of the solution of these systems is analyzed.

## C. Two Definitions for Stochastic Boundedness

As discussed in the introduction section, system (3) with nonvanishing PK does not possess any equilibrium point. Thus, stability analysis is meaningless and the boundedness property should be used to study the behavior of this system. In this section, two definitions for the stochastic boundedness of a stochastic process are presented.

**Definition 1** [12]: The stochastic process $X(t,\omega)$, defined on the probability space $(\Omega, \Im, \mathrm{P})$ and time set $\mathrm{T}$, is said to be $\ell^p$ *bounded* if there exists a real number $0 \leq b < \infty$, such that.

$$\sup_{t \in \mathrm{T}} E\left(\|X(t,\omega)\|_p\right) < b \tag{7}$$
$$\blacksquare$$

**Definition 2** [5]: The stochastic process $X(t,\omega)$, defined on the probability space $(\Omega, \Im, \mathrm{P})$ and time set $\mathrm{T}$, is said to be $p-$ *bounded in probability* if

$$\lim_{r \to \infty} \sup_{t \in \mathrm{T}} \mathrm{P}\left\{\|X(t,\omega)\|_p > r\right\} = 0 \tag{8}$$
$$\blacksquare$$

Although both of these definitions consider the boundedness issue of a stochastic process and coincide in many cases, they study different aspects of this phenomenon. In Definition 1, the expectation value of the process $p$ norm is considered. This property states that the trajectory does not stochastically leave a bound, containing the origin, in the sense of $p$ norm. The aim of the authors to propose this boundedness definition is to find the upper bound $b$ for stochastically perturbed systems in upcoming sections. Likewise, in Definition 2, the boundedness phenomenon is studied in the sense of probability distribution of the process. In the following sections, this definition is employed in order to analyze the boundedness of the state trajectory of the nonvanishing perturbed LTI systems and also to analyze the probabilistic distribution.

The aim of this paper is to define conditions on the system (3) such that the boundedness of its solutions, in the sense of the proposed definitions, is achieved. In the next section, it is shown that under some conditions on $A$ and $g(X,t)$, the strong solution is stochastically bounded, and the bounds are evaluated.

# III. BOUNDEDNESS OF NONVANISHING PERTURBED LTI SYSTEMS

As previously discussed, for stable systems with nonvanishing perturbation of the form (3), finding a bound on the state trajectory is of tremendous importance.

Let us define the following parameter:

$$c \triangleq \sup_{t \in T} \|g(0,t)\| \tag{9}$$

which is finite due to Assumption 1. Also, based on the triangular inequality and Lipchitz continuity condition of PK in Asumption 1, we have:

$$\|g(X,t)\| \leq c + \gamma \|X\| \tag{10}$$

In the following, a theorem is stated, which shows that under a condition on the Lyapunov complement and the Lipchitz constant $\gamma$ defined in Assumption 1.2, how the stochastic system's response $X(t,\omega)$ remains stochastically bounded, and the bounds are calculated.

**Theorem 1**: Consider a nonlinearly nonvanishing perturbation of an asymptotically stable LTI system in the form of (3), if the following inequality is satisfied:

$$\gamma^2 \bar{\lambda}_P < \underline{\lambda}_Q \tag{11}$$

where $P$ is the Lyapunov complement of $A$ in (1), the following statements hold:

A. Each state trajectory starting almost surely from $X_0$ approaches in $\ell^2$ norm to the following bound:

$$\lim_{t \to \infty} \{E(\|X(t)\|)\} \leq b = \frac{-c\gamma \bar{\lambda}_P - \sqrt{c^2 \bar{\lambda}_P \underline{\lambda}_Q}}{\gamma^2 \bar{\lambda}_P - \underline{\lambda}_Q} \tag{12}$$

B. Each state trajectory starting almost surely from $X_0$ is $\ell^2$ bounded with the following bound:

$$\sup_{t \in T} \{E(\|X(t)\|)\} \leq \max \left\{ E(\|X_0\|), b = \frac{-c\gamma \bar{\lambda}_P - \sqrt{c^2 \bar{\lambda}_P \underline{\lambda}_Q}}{\gamma^2 \bar{\lambda}_P - \underline{\lambda}_Q} \right\} \tag{13}$$

C. Each state trajectory starting almost surely from $X_0$ with $b < E(\|X_0\|)$ is $2-$ bounded in probability for any $b < \varepsilon$.

$$P\left(\sup_{t \in T} \|X\| > \varepsilon\right) \leq \frac{1}{\varepsilon} \underline{\lambda}_P \overline{\lambda}_P E(\|X_0\|) \tag{14}$$

∎

**Proof**:

A. Consider the following quadratic time-invariant Lyapunov function for system (3):

$$V(X) = \frac{1}{2} X^T P X \tag{15}$$

where $P$ is the Lyapunov compliment in (1). Note that $V$ is continuously twice differentiable in $X$, $V \in C^2(\mathbb{R}^n)$. Using the Ito formula [3, 11, 13], the stochastic increment of this Lyapunov function due to state dynamics (3) is:

$$dV(X) = \left(X^T P A X + \frac{1}{2} g^T P g\right) dt + X^T P g dW \tag{16}$$

Using the fact that the first term is scalar and $P$ is symmetric, (16) may be rewritten as follows due to (1):

$$\begin{aligned} dV(X) &= \left(\frac{1}{2} X^T (PA + A^T P) X + \frac{1}{2} g^T P g\right) dt + X^T P g dW \\ &= \left(\frac{-1}{2} X^T Q X + \frac{1}{2} g^T P g\right) dt + X^T P g dW \end{aligned} \tag{17}$$

Using sub-multiplicativity property and triangular norm inequality along with reordering the terms, we obtain:

$$dV(X) \leq \begin{pmatrix} \dfrac{-\underline{\lambda}_Q + \gamma^2 \overline{\lambda}_P}{2} \|X\|^2 \\ + c \gamma \overline{\lambda}_P \|X\| \\ + \dfrac{c^2 \overline{\lambda}_P}{2} \end{pmatrix} dt + X^T P g dW \tag{18}$$

where (10) has been used. Let us consider the expected value of $dV$ with respect to probability P defined as the following Lebesgue integral on the sample space $\Omega$ [14]:

$$E(dV) \triangleq \int_\Omega dV(X(\omega))dP(\omega) \leq \begin{pmatrix} \dfrac{-\lambda_Q + \gamma^2 \overline{\lambda}_P}{2} E(\|X\|^2) \\ +c\gamma\overline{\lambda}_P E(\|X\|) \\ +\dfrac{c^2 \overline{\lambda}_P}{2} \end{pmatrix} dt \quad (19)$$

where the fact that the standard Wiener process has zero mean is used.

Assume that the term $\dfrac{-\lambda_Q + \gamma^2 \overline{\lambda}_P}{2}$, which is the coefficient of $E(\|X\|^2)$, is negative. Also, note that the function $f(x) = x^2$ is a convex function. Now using the Jensen's inequality yields [13, 14]:

$$f\left(\int_\Omega \|X\| dP\right) \leq \int_\Omega f(\|X\|) dP$$
$$\Rightarrow E^2(\|X\|) \leq E(\|X\|^2) \quad (20)$$
$$\Rightarrow \dfrac{-\lambda_Q + \gamma^2 \overline{\lambda}_P}{2} E(\|X\|^2) \leq \dfrac{-\lambda_Q + \gamma^2 \overline{\lambda}_P}{2} E^2(\|X\|)$$

Consequently, we obtain:

$$E(dV) \leq \begin{pmatrix} \dfrac{-\lambda_Q + \gamma^2 \overline{\lambda}_P}{2} E^2(\|X\|) \\ +c\gamma\overline{\lambda}_P E(\|X\|) \\ +\dfrac{c^2 \overline{\lambda}_P}{2} \end{pmatrix} dt \triangleq L_V(X) dt \quad (21)$$

Consider the inequality $L_V(X) \leq 0$, obviously it is quadratic in $E(\|X(t)\|)$ with the discriminant:

$$\Delta = c^2 \overline{\lambda}_P \lambda_Q > 0 \quad (22)$$

which expresses that the equality $L_V(X) = 0$ has two distinct solutions in $E(\|X(t)\|)$ (considering all solutions including those that are not allowed). Using the fact that this quadratic form has downside curvature ($\dfrac{-\lambda_Q + \gamma^2 \overline{\lambda}_P}{2} < 0$), the equality $L_V(X) = 0$ has the following two solutions:

$$E(\|X\|)_1 = b = \frac{-c\gamma\bar{\lambda}_P - \sqrt{c^2\bar{\lambda}_P \underline{\lambda}_Q}}{\gamma^2\bar{\lambda}_P - \underline{\lambda}_Q} > 0$$

$$E(\|X\|)_2 = \frac{-c\gamma\bar{\lambda}_P + \sqrt{c^2\bar{\lambda}_P \underline{\lambda}_Q}}{\gamma^2\bar{\lambda}_P - \underline{\lambda}_Q} < 0 \tag{23}$$

Figure 1 shows an upper bound for $L_V(X)$ with respect to $E(\|X(t)\|)$:

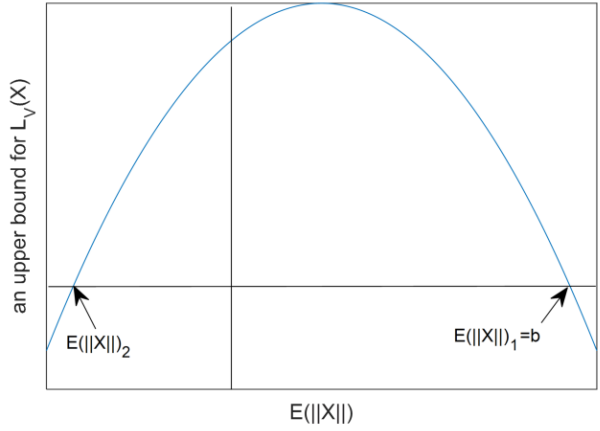

Figure 1- Upper bound for $L_V(X)$ with respect to $E(\|X(t)\|)$.

Now consider the following subsets:

$$S_o \triangleq \{X(t,\omega) \in \mathbb{R}^n \,|\, E(\|X(t)\|) > E(\|X\|)_1\}$$
$$S_i \triangleq \{X(t,\omega) \in \mathbb{R}^n \,|\, E(\|X(t)\|) \leq E(\|X\|)_1\} \tag{24}$$

Thus,

$$L_V(X) < 0, \forall X \in S_o \quad \text{a.s} \tag{25}$$

and:

$$E(dV(X)) < 0, \forall X \in S_o \quad \text{a.s} \tag{26}$$

Let us use the Dynkin's theorem [13] to analyze the time evolution of the Lyapunov value. If $t > s \in T$, and $t - s$ is sufficiently small, then Dynkin's formula yields:

$$E\left(V(t)-V(s)|\Im(W(s))\right)=\int_s^t E\left(dV(\tau)\right)d\tau<0 \quad {}^2 \qquad (27)$$
$$\Rightarrow E\left(V(t)|\Im(W(s))\right)<E\left(V(s)\right)$$

which shows that $V(t)$ is a supermartingale. Now assume that $S_o$ is an invariant set of the states, i.e., if $X(t_0)\in S_o$ then, $\forall t_1>t_0 : X(t_1)\in S_o$. On the other hand, the Lyapunov value is decreasing in time, and the fact that $\|X\|\leq \frac{V(X)}{\underline{\lambda}_p}$ contradicts with $X(t)\in S_o, t>t_0$ to be bounded from bellow in Euclidean norm. Thus, $S_o$ is not an invariant set of states. Regarding the existence and uniqueness of solution, the system approaches to its compact complement, $S_i$. Thus:

$$\lim_{t\to\infty}\{E(\|X(t)\|)\}\leq b = \frac{-c\gamma\bar{\lambda}_P - \sqrt{c^2\bar{\lambda}_P\underline{\lambda}_Q}}{\gamma^2\bar{\lambda}_P - \underline{\lambda}_Q} \qquad (28)$$

This argument results A.

B. On the other hand, for the trajectories starting from $S_o$, i.e., $E(\|X_0\|)>b$, the value of $E(\|X\|)$ decreases in time until it approaches $S_i$. Also, the trajectories starting from $S_i$, remain there. Therefore, an upper bound for $E(\|X\|)$ is achieved as the following:

$$\sup_{t\in T}\{E(\|X(t)\|)\}\leq \max\left\{E(\|X_0\|), b = \frac{-c\gamma\bar{\lambda}_P - \sqrt{c^2\bar{\lambda}_P\underline{\lambda}_Q}}{\gamma^2\bar{\lambda}_P - \underline{\lambda}_Q}\right\} \qquad (29)$$

which infers B.

C. Regarding the fact that $V(t)$ is supermartingale, if $X\in S_o$, which is inferred from (17), the supermartingale inequality [3] yields that for every $b<\varepsilon$:

$$P\left(\sup_{t\in T} V(X)>\varepsilon\right)\leq \frac{1}{\varepsilon}\left(E(V(0))+\max(0, E\lim_{t\to\infty}(-V(t)))\right) \qquad (30)$$

Thus,

$$P\left(\sup_{t\in T} V(X)>\varepsilon\right)\leq \frac{1}{\varepsilon}V(E(X(0)))\leq \frac{1}{\varepsilon}\bar{\lambda}_P E(\|X_0\|) \qquad (31)$$

and (32)

---

[2] $V(t)$ is abbreviation for $V(X(t,\omega))$

$$P\left(\sup_{t\in T}\|X\| > \varepsilon/\underline{\lambda}_P\right) \leq \frac{1}{\varepsilon}\bar{\lambda}_P E(\|X_0\|)$$

where the fact that $\underline{\lambda}_P\|X\| \leq V(X) \leq \bar{\lambda}_P\|X\|$ was used to derive the last two formulae. Now, let us define $\hat{\varepsilon} \triangleq \frac{\varepsilon}{\underline{\lambda}_P}$, we have:

$$P\left(\sup_{t\in T}\|X\| > \hat{\varepsilon}\right) \leq \frac{1}{\hat{\varepsilon}}\underline{\lambda}_P \bar{\lambda}_P E(\|X_0\|) \tag{33}$$

Due to the fact that $\underline{\lambda}_P \bar{\lambda}_P \|X_0\|$ is finite, in the limiting situation of $\hat{\varepsilon} \to \infty$, (14) is verified. ∎

By this theorem, the stochastic systems response under the proposed conditions on the system, is both $\ell^2$ bounded and $2-$ bounded in probability.

**Remark 1**: Consider the inequality $\frac{-\underline{\lambda}_Q + \gamma^2 \bar{\lambda}_P}{2} < 0$ in Theorem 1. We want to know if for a given $\gamma$, a pair of matrices $P$ and $Q$ can be found in a way that the condition $\frac{-\underline{\lambda}_Q + \gamma^2 \bar{\lambda}_P}{2} < 0$ holds. Using the fact that both $P$ and $Q$ are positive definite matrices, we consider the following two cases for $\gamma$. First consider the case $\gamma = 0$, which yields:

$$\gamma = 0 \Rightarrow \frac{-\underline{\lambda}_Q + \gamma^2 \bar{\lambda}_P}{2} = \frac{-\underline{\lambda}_Q}{2} < 0 \tag{34}$$

So, for $\gamma = 0$, every choice of $Q$ is suitable and acceptable and guarantees $\frac{-\underline{\lambda}_Q + \gamma^2 \bar{\lambda}_P}{2} < 0$. It is notable that this case models constant nonzero PKs i.e., $g = $ cte.

Second, consider the case $\gamma \neq 0$ for which we have:

$$\gamma \neq 0 \Rightarrow \frac{-\underline{\lambda}_Q + \gamma^2 \bar{\lambda}_P}{2} < 0 \Leftrightarrow \frac{\bar{\lambda}_P}{\underline{\lambda}_Q} < \frac{1}{\gamma^2} \tag{35}$$

Using the fact that the ratio $\frac{\bar{\lambda}_P}{\underline{\lambda}_Q}$ obtains its minimum value when $Q = I$. For $\frac{-\underline{\lambda}_Q + \gamma^2 \bar{\lambda}_P}{2} < 0$, it suffices to have the admissible pair of $P$ and $Q = I$, when

$$\gamma^2 \bar{\lambda}_P < 1 \tag{36}$$

Also, consider the bound $E(\|X\|)_1$ in $S_o$. Using the notation $\kappa \triangleq \frac{\bar{\lambda}_P}{\underline{\lambda}_Q}$, equation (23) can be rewritten as:

$$E(\|X\|)_1 = \frac{-c\gamma\kappa - \sqrt{c^2\kappa}}{\gamma^2\kappa - 1} \tag{37}$$

The partial ddifferivative of $E(\|X\|)_1$ with respect to $\kappa$ is:

$$\frac{\partial}{\partial \kappa} E(\|X\|)_1 = \frac{c\sqrt{\kappa}\left(\gamma + \kappa^{-\frac{1}{2}}\right)^2}{2(\gamma^2\kappa - 1)^2} > 0 \tag{38}$$

which intuitively infers that the bound will achieve its smallest value as $\kappa$ is minimized. Hence, choosing $Q = I$ gives the minimum value for the bound $b$ (denoted by $\underline{b}$). Thus, the system trajectory approaches the invariant set $S_i$ almost surely:

$$\lim_{t \to \infty} X(t) \in S_i = \left\{ X(\omega) \in \mathbb{R}^n \left| E(\|X\|) \leq \underline{b} = \frac{c\gamma\overline{\lambda}_{\hat{p}} + \sqrt{2c^2\overline{\lambda}_{\hat{p}}}}{1 - \gamma^2\overline{\lambda}_{\hat{p}}} \right. \right\} \tag{39}$$

Consequently, if Assumption 1 and (36) are satisfied, the system is $\ell^2$ bounded almost surely and the bound $\underline{b}$ is given by (39) and also it is $2-$ bounded in probability. ∎

## IV. AN ILLUSTRATIVE EXAMPLE

In order to show the effectiveness of the proposed method in finding a bound for the stochastic state trajectories of stable LTI systems in the presence of nonvanishing perturbation, consider the following system:

$$\dot{x} = -x + 0.25\cos(4x)dW \tag{40}$$

The intrinsic LTI system is stable and the PK has the following properties due to Assumption 1:

$$\begin{aligned} g(x,t) &= 0.25\cos(4x) \\ c &= 0.25 \\ \gamma &= 1 \end{aligned} \tag{41}$$

Also, the trivial Lyapunov equation has the unique solution $p = \frac{1}{2}$ and condition (11) is satisfied.

First, consider the trajectories starting from the origin. Theorem 1 predicts that the trajectories start from the origin do not leave the following bound:

$$\sup_{t \in T}\{E(\|X(t)\|)\} \leq b = 0.25 \tag{42}$$

Using Monte-Carlo simulation, 100 trajectories were generated. The expectation value $E(\|X(t)\|)$ and the predicted bound are plotted in Figure 2.

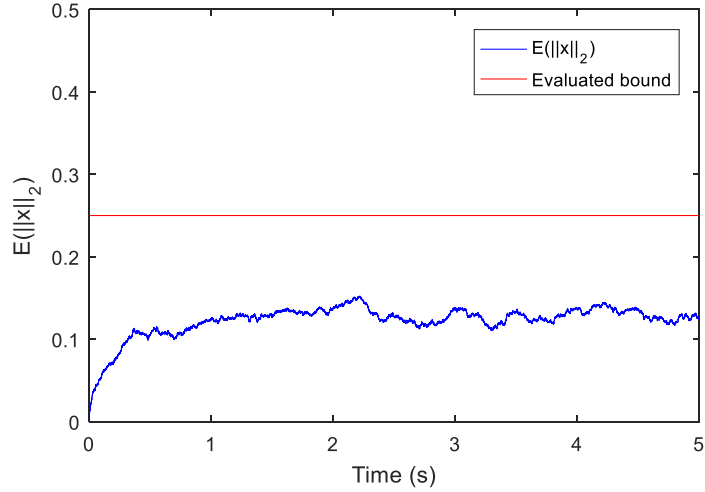

Figure 2- Expected value of state trajectories starting from the origin.

For the trajectories that do not start from the origin (trajectories start from $x(0) = 0.5$ in this simulation), since the intrinsic system is stable, the trajectories enter the prescribed bound given in (12) and (29). Also, the supremum does not exceed the bound in (13). We have used the Monte-Carlo simulation method and 100 trajectories have been used to evaluate the expected value, which is illustrated in Figure 3:

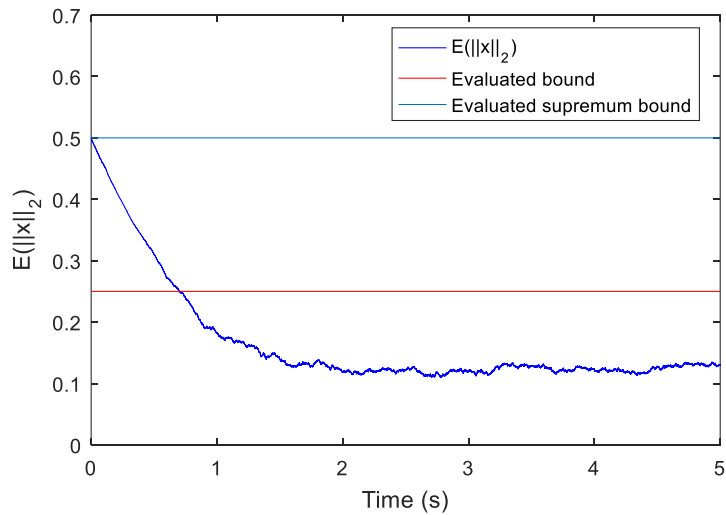

Figure 3- Expected value of state trajectories starting from $x(0) = 0.5$.

In order to illustrate part C of Theorem 1, 1000 simulations were performed, and the approximate probability of $\sup_{t \in \mathrm{T}} \|X\|$ as well as the evaluated bound in (14) are illustrated in Figure 4.

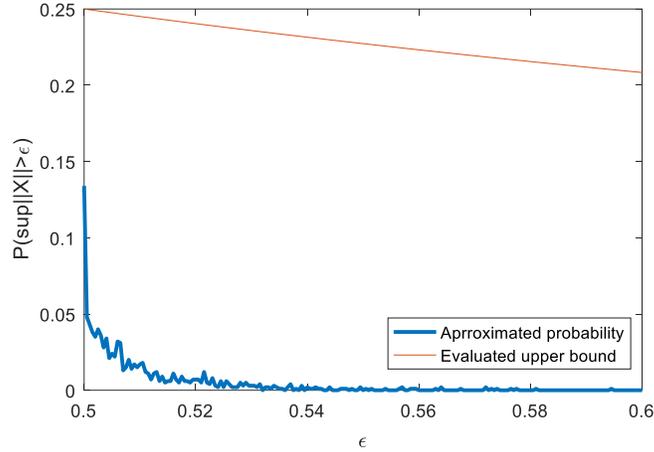

Figure 4-An upper bound for $\mathrm{P}\left(\sup_{t \in \mathrm{T}} \|X\| > \varepsilon\right)$ and its aproximated value obtained from Monte-Carlo simulation.

## V. CONCLUSION

In this paper, it was shown that if a stable LTI system is perturbed by a nonlinear stochastic perturbation, it may lose its stability, but under a mild condition on the Lipchitz constant of the perturbation kernel, it will remain stochastically bounded near its stable equilibrium point. This is a fairly common situation in experiments and the calculated bound helps designers to guarantee a level of tolerance for the state trajectories in the presence of Gaussian noise.